# Some results and conjectures about recurrence relations for certain sequences of binomial sums.

*Johann Cigler*

Fakultät für Mathematik
Universität Wien
A-1090 Wien, Nordbergstraße 15

Johann Cigler@univie.ac.at

**Abstract**

*In a previous paper [1] I have obtained some results about recurrences of certain binomial sums. The aim was to find an explanation of the remarkable identities*

$$\sum_{k \in \mathbb{Z}} (-1)^k \binom{n}{\left\lfloor \frac{n+5k+2}{2} \right\rfloor} = F_n \text{ and } \sum_{k \in \mathbb{Z}} (-1)^k \binom{n}{\left\lfloor \frac{n+5k}{2} \right\rfloor} = F_{n+1} \text{ for Fibonacci numbers } F_n \text{ by}$$

*putting these identities into a more general context. These identities have first been obtained by I. Schur [4] in his proof of the famous Rogers-Ramanujan identities. In the meantime computer experiments have led to several conjectures about the concrete form of these recurrences, some of which are proved in this paper.*

**1. Some general results**

We recall the following known fact (cf. [1]):

**Theorem 1**
Let $m \geq 2, i \geq 1$ be integers, $n \in \mathbb{N}$ and $\ell \in \mathbb{Z}$.
The sequences

$$a(n,i,\ell,m,z) = \sum_{k \in \mathbb{Z}} \binom{n}{\left\lfloor \frac{n+ik+\ell}{m} \right\rfloor} z^k \in \mathbb{Q}\left[z, z^{-1}\right] \tag{1}$$

satisfy linear recurrences of order $i$ with constant coefficients. More precisely there exist uniquely determined polynomials

$$p_1(n,m,x,s) = x^n + \sum_{j=1}^{\left\lfloor \frac{n}{m} \right\rfloor} c_1(n,m,j) x^{n-mj} s^j \tag{2}$$

and

$$p_k(n,m,x,s) = \sum_{j=1+\left\lfloor \frac{(k-1)n}{m} \right\rfloor}^{\left\lfloor \frac{kn}{m} \right\rfloor} c_k(n,m,j) x^{kn-mj} s^j, \tag{3}$$

with integer coefficients such that

$$\sum_{k=1}^{m-1} z^{k-1} p_k(i,m,E,1) a(n,i,\ell,m,z) = \left( \frac{1}{z} + (-1)^{m(i-1)} z^{m-1} \right) a(n,i,\ell,m,z), \tag{4}$$

where $E$ denotes the shift operator defined by $E^k f(n) = f(n+k)$.

### Sketch of the proof
We observe as in [3] that the recurrence relations for the binomial coefficients imply that

$$t(m,n,k) = \binom{n}{\left\lfloor \frac{n+k}{m} \right\rfloor}$$ satisfies the recurrence relation

$$t(m,n,k) = t(m,n-1,k-m+1) + t(m,n-1,k+1). \tag{5}$$

If we introduce the shift operator $K$ defined by $K^j f(k) = f(k-j)$, then (5) can be written in the form

$$t(m,n,k) = Et(m,n-1,k) = (K^{m-1} + K^{-1})t(m,n-1,k) = (K^{m-1} + K^{-1})^n t(m,0,k). \tag{6}$$

This implies

$$p_j(i,m,E,1)t(m,n,k) = p_j(i,m,K^{m-1}+K^{-1},1)t(m,n,k). \tag{7}$$

It turns out that in order to prove (4) it suffices to show

$$\sum_{j=1}^{m-1} x^{(j-1)n} p_j(i,m,x^{m-1}+\frac{1}{x},1) = \frac{1}{x^i} + (-1)^{m(i-1)} x^{(m-1)i}. \tag{8}$$

For then we get from (7)

$$\sum_{j=1}^{m-1} p_j(i,m,E,1)t(m,n,k-(j-1)i) = t(m,n,k+i) + (-1)^{m(i-1)} t(m,n,k-(m-1)i) \tag{9}$$

for all $n \in \mathbb{N}$.

This means that

$$\sum_{j=1}^{m-1} p_j(i,m,E,1)\left(\left\lfloor\frac{n+k+\ell-(j-1)i}{m}\right\rfloor\right) = \left(\left\lfloor\frac{n+k+i+\ell}{m}\right\rfloor\right) + (-1)^{m(i-1)}\left(\left\lfloor\frac{n+k-(m-1)i+\ell}{m}\right\rfloor\right).$$

If we replace $k$ by $ik$, multiply each side by $z^k$ and sum over all $k \in \mathbb{Z}$ we get

$$\sum_{j=1}^{m-1} z^{j-1} p_j(i,m,E,1) a(n,i,\ell,m,z) = \left(\frac{1}{z} + (-1)^{m(i-1)} z^{m-1}\right) a(n,i,\ell,m,z),$$

Identity (8) is equivalent with

$$\sum_{k=1}^{m-1} p_k(n,m,x+1,x) = 1 + (-1)^{m(n-1)} x^n. \tag{10}$$

This follows from

$$\sum_{k=1}^{m-1} x^{kn} p_k(n,m,x^{m-1}+\frac{1}{x},1) = \sum_{k=1}^{m-1}\sum_j c_k(n,m,j) x^{kn}\left(x^{m-1}+\frac{1}{x}\right)^{kn-mj}$$

$$= \sum_{k=1}^{m-1}\sum_j c_k(n,m,j) x^{mj}\left(x^m+1\right)^{kn-mj} = \sum_{k=1}^{m-1} p_k(n,m,x^m+1,x^m) = 1 + (-1)^{m(n-1)} x^{mn}.$$

For $n=0$ we set $p_k(0,m,x,s) = (-1)^{k-1}\binom{m}{k}$. Then (10) remains true for $n=0$.

It is easy to verify that the polynomials $p_k(n,m,x,s)$ of the form (3) which satisfy (10) are uniquely determined.

Let now $b(n,m,j) = (-1)^{m(n-1)-1} c_k(n,m,j)$ for $\left\lfloor\frac{(k-1)n}{m}\right\rfloor < j \leq \frac{kn}{m}$.

Then (10) may be reformulated as

$$(-1)^{m(n-1)} x^n + (-1)^{m(n-1)} \sum_{k=1}^{m-1}\sum_{j=1+\left\lfloor\frac{(k-1)n}{m}\right\rfloor}^{\left\lfloor\frac{kn}{m}\right\rfloor} b(n,m,j)(x+1)^{kn-mj} x^j = (1+x)^n - 1. \tag{11}$$

If we consider the homomorphism which sends $1+x$ to a primitive root of unity $\zeta_n$ of order $n$, then it has been shown in [1] that (11) implies

$$\sum_{j=1}^{n} b(n,m,j)x^j = \prod_{k=1}^{n}\left(x - \zeta_n^{-mk}(\zeta_n^k - 1)\right). \qquad (12)$$

Let $s_j$ be the $j$'s elementary symmetric function of $\zeta_n^{-mk}(\zeta_n^k - 1)$, $1 \le j < n$, and $\pi_j$ the $j$'s power sum of those numbers. Let $k(\mod n)$ be the remainder modulo $n$ with $0 \le k(\mod n) < n$.

Then $\pi_j = \sum_{i=0}^{j}(-1)^{j-i}\binom{j}{i}\sum_{k=1}^{n}\zeta_n^{ki-mkj} = (-1)^{j-mj(\mod m)} n \binom{j}{mj(\mod n)}.$

Newton's formula gives

$\sum_{i=0}^{j-1}(-1)^i \pi_{j-i} s_i = (-1)^{j-1} j s_j.$

Since $s_j = (-1)^j b_{n,m,n-j}$ we get

**Theorem 2**

*Let* $d(n,m,j) = (-1)^{(j-mj(\mod n))}\binom{j}{mj(\mod n)} n$ *and let* $b(n,m,n-j)$ *be defined by*

$b(n,m,n-j) = -\frac{1}{j}\sum_{i=0}^{j-1} d(n,m,j-i) b(n,m,n-i)$ *with* $b(n,m,n) = 1.$

*Then for* $n \ge 1$ *we have*

$$p_1(n,m,x,s) = x^n + \sum_j c_1(n,m,j)x^{n-mj}s^j = x^n - (-1)^{m(n-1)}\sum_{j=1}^{\left\lfloor \frac{n}{m} \right\rfloor} b(n,m,j)x^{n-mj}s^j \qquad (13)$$

*and*

$$p_k(n,m,x,s) = \sum_j c_k(n,m,j)x^{kn-mj}s^j = -(-1)^{m(n-1)}\sum_{j=1+\left\lfloor \frac{(k-1)n}{m} \right\rfloor}^{\left\lfloor \frac{kn}{m} \right\rfloor} b(n,m,j)x^{kn-mj}s^j. \qquad (14)$$

## 2. Concrete evaluations

Now we want to obtain more concrete information about the polynomials $p_k(n,m,x,s)$.
First I derive some known results (cf. [1]) from the present point of view.

### 2.1. The case $m=2$

For $m=2$ formula (10) tells us that there are polynomials $p_1(n,2,x,s)$ of the form $\sum_j c_j x^{n-2j} s^j$ which satisfy

$$p_1(n,2,x+1,x) = x^n + 1. \tag{15}$$

It is easy to see that the sequence $\left(p_1(n,2,x,1)\right)_{n\geq 1}$ begins with
$x, x^2-2, x^3-3x, x^4-4x^2+2, x^5-5x^3+5x, x^6-6x^4+9x^2-2,\cdots$.

Everyone familiar with the Lucas polynomials $L_n(x,s) = \sum_{j=0}^{\lfloor\frac{n}{2}\rfloor} \binom{n-j}{j} \frac{n}{n-j} x^{n-2j} s^j$ will immediately guess that

$$p_1(n,2,x,1) = \sum_{j=0}^{\lfloor\frac{n}{2}\rfloor} (-1)^j \binom{n-j}{j} \frac{n}{n-j} x^{n-2j} \tag{16}$$

and thus
$p_1(n,2,x,s) = L_n(x,-s)$ holds.
It is well known that the Lucas polynomials satisfy the recurrence
$L_n(x,s) = xL_{n-1}(x,s) + sL_{n-2}(x,s)$ with initial values $L_0(x,s) = 2$ and $L_1(x,s) = x$.

Therefore it remains only to show that $p_1(n,2,x,s) = L_n(x,-s)$ satisfies (15), i.e.

$$L_n(x+1,-x) = x^n + 1. \tag{17}$$

This is of course a well-known result. It may be proved in the following way:
The polynomials $L_n(x+1,-x)$ satisfy the recurrence
$L_{n+2}(x+1,-x) - (x+1)L_{n+1}(x+1,-x) + xL_n(x+1,-x) = \left(E^2 - (x+1)E + x\right) L_n(x+1,-x) = 0$.
If $p(E)f(n) = 0$ we call $p(z)$ the characteristic polynomial of the sequence $f(n)$. Thus the characteristic polynomial of the sequence $L_n(x+1,-x)$ is $z^2 - (x+1)z + x = (z-1)(z-x)$.
This implies that the constant sequence 1 and the sequence $(x^n)$ satisfy the same recurrence as $\left(L_n(x+1,-x)\right)$.
Therefore we get
$\left(E^2 - (x+1)E + x\right)\left(L_n(x+1,-x) - 1 - x^n\right) = 0$
with initial values $L_0(x+1,-x) - 1 - x^0 = 2-1-1 = 0$ and
$L_1(x+1,-x) - 1 - x = (x+1) - 1 - x = 0$.

This implies that $L_n(x+1,-x) - 1 - x^n = 0$ for all $n \in \mathbb{N}$ as asserted.

Another way to derive this result is by considering the generating function
$$\sum_{n \geq 0} p_1(n, 2, x+1, x) z^n = \sum_{n \geq 0} (1 + x^n) z^n = \frac{1}{1-z} + \frac{1}{1-xz} = \frac{2-(x+1)z}{1-(x+1)z + xz^2}.$$ Since the polynomials
are uniquely determined we get the generating function of the Lucas polynomials
$$\sum_{n \geq 0} p_1(n, 2, x, s) z^n = \frac{2 - xz}{1 - xz + sz^2}.$$

We have thus shown that
$$L_i(E, -1) a(n, i, \ell, 2, z) = \left(z + \frac{1}{z}\right) a(n, i, \ell, 2, z) \tag{18}$$
holds.

E.g. for $i = 5$ we get $L_5(E, -1) = E^5 - 5E^3 + 5E$. Therefore
$(E^5 - 5E^3 + 5E + 2) a(n, 5, \ell, 2, -1) = 0$. Now we have
$(E^5 - 5E^3 + 5E + 2) = (E+2)(E^2 - E - 1)^2$. The Fibonacci numbers satisfy the same
recurrence. Considering the initial values we get $\sum_{k \in \mathbb{Z}} (-1)^k \left( \left\lfloor \dfrac{n}{\frac{n+5k+2}{2}} \right\rfloor \right) = F_n$ and

$\sum_{k \in \mathbb{Z}} (-1)^k \left( \left\lfloor \dfrac{n}{\frac{n+5k}{2}} \right\rfloor \right) = F_{n+1}$.

In [3] we have also obtained the generating function
$$\sum_{n \geq 0} a(n, i, 0, 2, z) x^n = \frac{\dfrac{x^{i-1}}{z} + F_i(1, -x^2) + xF_{i-1}(1, -x^2)}{L_i(1, -x^2) - x^i \left(z + \dfrac{1}{z}\right)}, \tag{19}$$

where $F_n(x, s) = \sum_{k=0}^{n-1} \binom{n-1-k}{k} s^k x^{n-2k-1}$ are the Fibonacci polynomials.

As shown in [1] and [3] for $z = \pm 1$ there are always simpler recurrences. The formulas are different depending whether $i$ is odd or even. For example
$$\sum_{n \geq 0} a(n, 2m+1, 0, 2, -1) x^n = \frac{F_{2m+1}(1, -x^2) + xF_{2m}(1, -x^2) - x^{2m}}{L_{2m+1}(1, -x^2) + 2x^{2m+1}} = \frac{F_m(1, -x^2)}{F_{m+1}(1, -x^2) - xF_m(1, -x^2)}.$$

**Remark**

It should be noted that $a(n, 2m+1, 0, 2, -1)$ has a simple combinatorial interpretation (cf. [2]). It is the number of the set $A(n,m)$ of all lattice paths in $\mathbb{R}^2$ of length $n$ which are contained in the strip $-m-1 < y < m$, start at the origin and consist of $\left\lfloor \frac{n}{2} \right\rfloor$ northeast steps $(1,1)$ and $\left\lfloor \frac{n+1}{2} \right\rfloor$ southeast steps $(1,-1)$. Define a peak as a vertex preceded by a northeast step and followed by a southeast step, and a valley as a vertex preceded by a southeast step and followed by a northeast step. The height of a vertex is its $y$–coordinate. The peaks with height at least 1 and the valleys with height at most $-2$ are called extremal points. Define the weight of a path $v$ by $w(v,t) = t^{d(v)}$ where $d(v)$ is the number of extremal points of the path $v$. The weight of a set of paths is the sum of the weights of all paths of the set.

It has been shown in [2] that if we set $\binom{n}{k} = 0$ for $n < 0$ then the weight of the set $A(n,m)$ is

$$w(A(n,m),t) = \sum_{j \in \mathbb{Z}} (-1)^j \sum_{\ell \geq |j|} \binom{\left\lfloor \frac{n+(2m-3)j}{2} \right\rfloor}{\ell - j} \binom{\left\lfloor \frac{n+1-(2m-3)j}{2} \right\rfloor}{\ell + j} t^\ell.$$

For $t = 1$ this reduces by Vandermonde's formula to $\sum_{j \in \mathbb{Z}} (-1)^j \binom{n}{\left\lfloor \frac{n+(2m+1)j}{2} \right\rfloor}$.

For example $w(A(n,2),t) = \sum_{k=0}^{\left\lfloor \frac{n}{2} \right\rfloor} \binom{n-2k}{k} t^k = F_{n+1}(1,t)$.

As shown in [2] with other methods for $m \geq 2$ the polynomial $w(A(n,m),t)$ satisfies the recurrence relation

$$\left( F_m(E^2 + 1 - t, -E^2) - (1+E) F_{m-1}(E^2 + 1 - t, -E^2) + E F_{m-2}(E^2 + 1 - t, -E^2) \right) w(A(n,m),t) = 0$$

of order $2m-2$ which for $t = 1$ reduces to the recurrence relation

$$\left( F_{m+1}(E,-1) - F_m(E,-1) \right) w(A(n,m),1) = \left( F_{m+1}(E,-1) - F_m(E,-1) \right) a(n, 2m+1, 0, 2, -1) = 0$$

of order $m$.

It should also be noted that $\left( a(n,5,0,2,-1) \right)_{n \geq 0}$, $\left( a(n,7,0,2,-1) \right)_{n \geq 0}$ and $\left( a(n,9,0,2,-1) \right)_{n \geq 0}$ are the sequences A000045, A028495 and A061551 of The On-Line Encyclopedia of Integer Sequences.

## 2.2. The case $m = 3$.

In [1] I have obtained explicit formulas for the recurrence relations of the sums

$$a(n,i,\ell,3,z) = \sum_{k\in\mathbb{Z}} \left(\left\lfloor \frac{n+ik+\ell}{3} \right\rfloor\right) z^k.$$

From Theorem 1 we know that there are polynomials $v_i(x)$ and $w_i(x)$ such that

$$\left(v_i(E) - zw_i(E) - \left(\frac{1}{z} + (-1)^{i-1} z^2\right)\right) a(n,i,\ell,z) = 0.$$

Computer experiments led to the table

|         | 1 | 2     | 3         | 4          | 5           | 6              | 7                  |
|---------|---|-------|-----------|------------|-------------|----------------|--------------------|
| $v_i(x)$ | $x$ | $x^2$ | $x^3 - 3$ | $x^4 - 4x$ | $x^5 - 5x^2$ | $x^6 - 6x^3 + 3$ | $x^7 - 7x^4 + 7x$ |
| $w_i(x)$ | 0 | $2x$  | 3         | $2x^2$     | $5x$        | $3 + 2x^3$     | $7x^2$             |

An inspection of these polynomials led to the conjecture that they satisfy the recurrences

$$v_n(x) = xv_{n-1}(x) - v_{n-3}(x) \tag{20}$$

with initial values

$$v_0(x) = 3, v_1(x) = x, v_2(x) = x^2 \tag{21}$$

and

$$w_n(x) = xw_{n-2}(x) + w_{n-3}(x) \tag{22}$$

with initial values

$$w_0(x) = 3, w_1(x) = 0, w_2(x) = 2x. \tag{23}$$

It is not difficult to determine these polynomials explicitly: For $n > 0$ we have

$$v_n(x) = \sum_{3j \leq n} (-1)^j \binom{n-2j}{j} \frac{n}{n-2j} x^{n-3j} \tag{24}$$

and

$$w_n(x) = \sum_{3j \leq 2n} \binom{n-j}{2n-3j} \frac{n}{n-j} x^{2n-3j}. \tag{25}$$

For the recurrences are easily verified and the initial values coincide for $n = 1, 2, 3$.

Thus we have experimental evidence for

**Theorem 3**

*For $i \geq 1$ the sequence*

$$a(n,i,\ell,3,z) := \sum_{k \in \mathbb{Z}} \left( \left\lfloor \frac{n+ik+\ell}{3} \right\rfloor \right) z^k \tag{26}$$

*satisfies the recurrence*

$$(v_i(E) - zw_i(E))a(n,i,\ell,z) - (\frac{1}{z} + (-1)^{i-1}z^2)a(n,i,\ell,z) = 0 \tag{27}$$

*for all $n \in \mathbb{N}$.*

**Remark**

In [1] unfortunately some typos corrupted the formulation of this theorem.

With the notations introduced above we have only to verify that the polynomials $p_i(n,3,x,s)$ are given by $v_n(x,s)$ and $w_n(x,s)$, i.e.

$$p_1(n,3,x,s) = v_n(x,s) = \sum_{3j \leq n} (-1)^j \binom{n-2j}{j} \frac{n}{n-2j} x^{n-3j} s^j \tag{28}$$

and

$$-p_2(n,3,x,s) = w_n(x,s) = \sum_{3j \leq 2n} \binom{n-j}{2n-3j} \frac{n}{n-j} x^{2n-3j} s^j. \tag{29}$$

Note that $v_n(x) = v_n(x,1)$ and $w_n(x) = w_n(x,1)$.

With other words we have to show that

$$v_n(1+x,x) - w_n(1+x,x) - 1 + (-x)^n = 0. \tag{30}$$

These polynomials satisfy the recurrences
$v_n(x,s) = xv_{n-1}(x,s) - sv_{n-3}(x)$ and $w_n(x,s) = xsw_{n-2}(x) + s^2 w_{n-3}(x)$ for all $n \in \mathbb{N}$ if we set $v_0(x,s) = w_0(x,s) = 3$.
The sequences $v_n(1+x,x)$ and the constant sequence 1 satisfy the recurrence
$(E^3 - (1+x)E^2 + x)v_n(1+x,x) = 0$ and $(E^3 - (1+x)E^2 + x)1 = 0$ and the sequences

$w_n(1+x, x)$ and $(-x)^n$ satisfy $(E^3 - (1+x)xE - x^2)w_n(1+x, x) = 0$ and
$(E^3 - (1+x)xE - x^2)(-x)^n = 0$.

Therefore $f_n(x) = v_n(1+x, x) - w_n(1+x, x) - 1 + (-x)^n$ satisfies the recurrence
$(E^3 - (1+x)E^2 + x)(E^3 - (1+x)xE - x^2)f_n(x) = 0$ of order 6. Since the first 6 initial values are
0 we get $f_n(x) = 0$ for all $n \in \mathbb{N}$.

It is again instructive to consider the generating functions

$$\sum_{n \geq 0} p_1(n, 3, x, s) z^n = \frac{3 - 2xz}{1 - xz + sz^3} \text{ and } \sum_{n \geq 0} p_2(n, 3, x, s) z^n = -\frac{3 - xsz^2}{1 - xsz^2 - s^2 z^3}.$$

This gives $\sum_{n \geq 0} p_1(n, 3, x+1, x) z^n = \frac{3 - 2(x+1)z}{1 - (x+1)z + xz^3} = \frac{1}{1-z} + \frac{2 - xz}{1 - xz - xz^2}$

and $\sum_{n \geq 0} p_2(n, 3, x+1, x) z^n = -\frac{3 - (x+1)xz^2}{1 - (x+1)xz^2 - x^2 z^3} = \frac{-1}{1+xz} - \frac{2 - xz}{1 - xz - xz^2}.$

These again imply

$$\sum_{n \geq 0} (p_1(n, 3, x+1, x) + p_2(n, 3, x+1, x)) z^n = \frac{1}{1-z} - \frac{1}{1+xz} = \sum_{n \geq 0} (1 + (-1)^{n-1} x^n) z^n.$$

In addition to Theorem 3 we compute the generating function

$$\sum_{n \geq 0} a(n, i, 0, 3, z) x^n = \frac{b_i(x) + c_i(x)z + (-1)^{i-1} x^{i-1} z^2}{v_i(1, x^3) - \frac{z}{x^i} w_i(1, x^3) - x^i \left(\frac{1}{z} + (-1)^{i-1} z^2\right)}. \quad (31)$$

Computer experiments suggest that $b_i(x)$ satisfies the same recurrence as $v_i(1, x^3)$, i.e.
$b_i(x) = b_{i-1}(x) - x^3 b_{i-3}(x)$.
The initial values are $b_1(x) = 1, b_2(x) = 1 + x, b_3(x) = 1 + x + x^2$.
The generating function of these polynomials is
$$\sum_{i \geq 0} b_{i+1}(x) t^i = \frac{1 + xt + x^2 t^2}{1 - t + x^3 t^3}.$$

Let $\frac{1}{1 - t + x^3 t^3} = \sum_{n \geq 0} r_n(x) t^n$. Then it is easily verified that $r_n(x) = \sum_{k=0}^{\lfloor n/2 \rfloor} \binom{n - 2k}{k} x^{3k}$. The first
values of the sequence $(r_n(x))_{n \geq 0}$ are $1, 1, 1, 1 - x^3, 1 - 2x^3, 1 - 3x^3, 1 - 4x^3 + x^6, \cdots$.
This gives $b_{i+1}(x) = r_i(x) + x r_{i-1}(x) + x^2 r_{i-2}(x)$.

The polynomials $c_i(x)$ satisfy the same recurrence as $\dfrac{w_i(1,x^3)}{x^i}$, i.e.

$c_i(x) = xc_{i-2}(x) + x^3 c_{i-3}(x)$.

The initial values are $c_1(x) = 1, c_2(x) = 1-x, c_3(x) = x + 2x^2$.

The generating function of these polynomials is

$$\sum_{i\geq 0} c_{i+1}(x) t^i = \frac{1+(1-x)t+2x^2 t^2}{1-xt^2-x^3 t^3}.$$

Let $\dfrac{1}{1-xt^2-x^3 t^3} = \sum_{n\geq 0} s_n(x) t^n$. Then $s_n(x) = \sum_{j=0}^{\left\lfloor \frac{n+1}{3} \right\rfloor} \binom{\left\lfloor \frac{n+1}{3} \right\rfloor + j}{3j - \varepsilon(n)} x^{n-3j+\varepsilon(n)}$,

where $\varepsilon(n) \equiv n \pmod 3$ with $\varepsilon(n) \in \{-1, 0, 1\}$. This can easily be verified by considering each of the cases $s_{3n}(x), s_{3n-1}(x), s_{3n+1}(x)$ separately. The first values of the sequence $(s_n(x))_{n\geq 0}$ are
$1, 0, x, x^3, x^2, 2x^4, x^3 + x^6, 3x^5, x^4 + 3x^7, \cdots$.

Another representation is $s_n(x) = \sum_{k=0}^{n} \binom{\frac{n+k}{3}}{k} x^{n-k} [n+k \equiv 0 \pmod 3]$.

Finally we get $c_{i+1}(x) = s_i(x) + (1-x) s_{i-1}(x) + 2x^2 s_{i-2}(x)$.

In order to prove (31) we note first that

$$\left( v_i(1,x^3) - \frac{z}{x^i} w_i(1,x^3) - x^i \left( \frac{1}{z} + (-1)^{i-1} z^2 \right) \right) \sum_{n\geq 0} a(n,i,0,3,z) x^n$$

is a polynomial $d_i(x,z)$ in $x$ of degree $< i$. For a polynomial $p(x) = \sum_k p_k x^k$ let

$P_i(p(x)) = \sum_{k=0}^{i-1} p_k x^k$.

Therefore

$$d_i(x,z) = P_i\left( \left( v_i(1,x^3) - \frac{z}{x^i} w_i(1,x^3) \right) \sum_{n\geq 0} a(n,i,0,3,z) x^n \right).$$

Since no factor contains negative powers of $z$ we see that the sum of the first $i$ terms of

$$P_i\left( v_i(1,x^3) \sum_{n=0}^{i-1} a(n,i,0,3,0) x^n \right) = d_i(x,0).$$

Let $a_i(x,z) = \sum_{n=0}^{i-1} a(n,i,0,3,z) x^n$.

Then

$$a_i(x, z) = \sum_{n=0}^{i-1}\binom{n}{\left\lfloor\frac{n}{3}\right\rfloor} x^n + z \sum_{n=\left\lfloor\frac{i-2}{2}\right\rfloor}^{i-1}\binom{n}{\left\lfloor\frac{n+i}{3}\right\rfloor} x^n + z^2 x^{i-1}. \tag{32}$$

This implies
$$d_i(x,0) - d_{i-1}(x,0) + x^3 d_{i-3}(x,0) = P_i\left(v_{i-1}(1, x^3)\left(a_i(x,0) - a_{i-1}(x,0)\right) + x^3 v_{i-3}(1, x^3)\left(a_i(x,0) - a_{i-3}(x,0)\right)\right) = 0.$$

Since the initial values coincide we get $d_i(x,0) = b_i(x)$.

Next we determine the coefficient of $z^2$ of $d_i(x, z)$. From (32) we get
$$[z^2]\left(v_i(1, x^3)\sum_{n=0}^{i-1} a(n, i, 0, 3, 0)x^n\right) = x^{i-1}.$$

On the other hand we have
$$-\frac{z}{x^{2i}} w_{2i}(1, x^3) z \sum_{n=\left\lfloor\frac{2i-2}{2}\right\rfloor}^{2i-1}\binom{n}{\left\lfloor\frac{n+2i}{3}\right\rfloor} x^n = -z^2\left(2x^i + \cdots\right)\left(x^{i-1} + \cdots\right) = -2x^{2i-1}z^2 + \cdots.$$

This implies $\left[z^2\right] d_{2i}(x, z) = -x^{2i-1}$.

For odd indices we get
$$-\frac{z}{x^{2i+1}} w_{2i+1}(1, x^3) z \sum_{n=\left\lfloor\frac{2i+1-2}{2}\right\rfloor}^{2i}\binom{n}{\left\lfloor\frac{n+2i+1}{3}\right\rfloor} x^n = -z^2 x^{i+2}(\cdots)\left(x^{i-1} + \cdots\right) = -x^{2i+1} z^2(\cdots).$$

This implies $\left[z^2\right] d_{2i+1}(x, z) = x^{2i-1}$.

As above we verify that

$$p_i(x) = P_i\left(\frac{1}{x^i} w_i(1, x^3)\sum_{n=0}^{i-1}\binom{n}{\left\lfloor\frac{n}{3}\right\rfloor} x^n\right)$$

satisfies the recurrence $p_i(x) - x p_{i-2}(x) - x^3 p_{i-3}(x) = 0$, i.e the same recurrence as $\frac{1}{x^i} w_i(1, x^3)$. In the same way we get that $q_i(x) = P_i\left(v_i(1, x^3)\sum_{n=\left\lfloor\frac{i-2}{2}\right\rfloor}^{i-1}\binom{n}{\left\lfloor\frac{n+i}{3}\right\rfloor} x^n\right)$ satisfies the same recurrence.

Writing $d_i(x, z) = b_i(x) + c_i(x) z + (-1)^{i-1} x^{i-1} z^2$ we see therefore that $c_i(x) - x c_{i-2}(x) - x^3 c_{i-3}(x) = 0$, from which equation (31) follows.

## 2.3. Some general observations.

As a special case of Theorem 1 we see that

$$a(n,1,0,m,z) = \sum_{k\in\mathbb{Z}} \left( \left\lfloor \frac{n+k}{m} \right\rfloor \right) z^k = \frac{1-z^m}{1-z}\left(\frac{1+z^m}{z}\right)^n.$$

From

$$p_1(n,m,x+1,x) = x^n + \sum_{j=1}^{\lfloor\frac{n}{m}\rfloor} c_1(n,m,j)(x+1)^{n-mj} x^j \text{ and the fact that}$$

$$p_k(n,m,x+1,x) = x^{1+\lfloor\frac{(k-1)n}{m}\rfloor}(\cdots) \text{ for } k \geq 2 \text{ we see that (10) implies that all coefficients}$$

$[x^j]p_1(n,m,x+1,x), 1 \leq j \leq \frac{n}{m}$, must vanish, i.e.

$P_{i+\lfloor\frac{n}{m}\rfloor}(p_1(n,m,x+1,x)) = 1$. By these conditions $p_1(n,m,x,s)$ is uniquely determined.

This observation implies that

$$p_1(n,m,x,s) = \sum_{j=0}^{\lfloor\frac{n}{m}\rfloor} (-1)^j \binom{n-(m-1)j}{j} \frac{n}{n-(m-1)j} x^{n-mj} s^j. \tag{33}$$

For

$$p_1(n,m,x+1,x) = \sum_{j=0}^{\lfloor\frac{n}{m}\rfloor} (-1)^j \binom{n-(m-1)j}{j} \frac{n}{n-(m-1)j} (x+1)^{n-mj} x^j$$

and therefore the coefficient of $x^r$ is given by

$$\sum_{j=0}^{r} (-1)^j \binom{n-(m-1)j}{j} \frac{n}{n-(m-1)j} \binom{n-mj}{r-j} = \frac{n}{r!}\sum_{j=0}^{r}(-1)^j\binom{r}{j}(n-(m-1)j-1)\cdots(n-(m-1)j-r+1).$$

for all $r$ such that $r \leq \left\lfloor\frac{n}{m}\right\rfloor$.

Let $q(n) = (n-1)_{r-1} = (n-1)(n-2)\cdots(n-r+1)$. This is a polynomial in $n$ of degree $r-1$.

With this notation the coefficient of $x^r$ is $\frac{n}{r!}(1-E^{m-1})^r q(n)$. Since

$(1-E^{m-1})n^k = n^k - (n+m-1)^k = k(m-1)n^{k-1} + \cdots$ the operator $1-E^{m-1}$ decreases the degree

of a polynomial by 1. Therefore $\frac{n}{r!}(1-E^{m-1})^r q(n) = 0$ and the assertion is proved.

It is easily verified that (33) satisfies the recurrence relation

$$\left(E^m - xE^{m-1} - s\right) p_1(n, m, x, s) = 0.$$

Now we observe that

$$p_{m-1}(n, m, x, s) = \sum_{j=0}^{\left\lfloor \frac{n}{m} \right\rfloor} \binom{n-j}{(m-1)n - mj} \frac{n}{n-j} x^{(m-1)n - mj} s^j.$$

We know that

$$d(n, m, j) = (-1)^{(j - mj \pmod{n})} \binom{j}{mj \pmod{n}} n.$$

This implies that $d(n, m, 0) = n$ and $d(n, m, j) = 0$ for $j \leq \frac{n}{m}$. For in this case we have

$mj \pmod{m} = mj$ and therefore $\binom{j}{mj} = 0.$

Thus also

$$b(n, m, n - j) = -\frac{1}{j} \sum_{i=0}^{j-1} d(n, m, j - i) b(n, m, n - i) = 0.$$

Let now $\frac{n}{m} < j \leq \frac{2n}{m}.$

From

$$-jb(n, m, n - j) = \sum_{i=0}^{j-1} d(n, m, j - i) b(n, m, n - i)$$

$$= d(n, m, j) b(n, m, n) + \sum_{i=1}^{\left\lfloor \frac{m}{n} \right\rfloor} d(n, m, j - i) b(n, m, n - i) + \sum_{i=1+\left\lfloor \frac{m}{n} \right\rfloor}^{j-1} d(n, m, j - i) b(n, m, n - i)$$

we deduce that $-jb(n, m, n - j) = d(n, m, j) = (-1)^{(j - mj \pmod{n})} \binom{j}{mj \pmod{n}} n,$

because the first sum vanishes since $b(n, m, n - i) = 0$ and the second sum vanishes because $d(n, m, j - i) = 0.$

Therefore for $\frac{(m-2)n}{m} < j \leq \frac{(m-1)n}{m}$ we get

$$b(n, m, j) = -(-1)^{(n - j + m(n-j) \pmod{n})} \binom{n-j}{m(n-j) \pmod{n}} \frac{n}{n-j} = -(-1)^{n + j + (m-1)n - mj} \binom{n-j}{(m-1)n - mj} \frac{n}{n-j}.$$

Therefore we get

$$p_{m-1}(n,m,x,s) = (-1)^m \sum_{j=1+\left\lfloor\frac{(m-2)n}{m}\right\rfloor}^{\left\lfloor\frac{(m-1)n}{m}\right\rfloor} \binom{n-j}{(m-1)n-mj}\frac{n}{n-j}x^{(m-1)n-mj}((-1)^{m-1}s)^j. \quad (34)$$

In order to give a simple description of these results we introduce polynomials

$$v_n(m,x,s) = \sum_{j\in\mathbb{N}} (-1)^j \binom{n-(m-1)j}{j}\frac{n}{n-(m-1)j}x^{n-mj}s^j \quad (35)$$

and

$$w_n(m,x,s) = \sum_{j\in\mathbb{N}} \binom{n-j}{(m-1)n-mj}\frac{n}{n-j}x^{(m-1)n-mj}s^j, \quad (36)$$

where as always in this paper $\binom{n}{k}=0$ for $n<0$.

These polynomials satisfy the recurrence relations

$$v_{n+m}(m,x,s) = xv_{n+m-1}(m,x,s) - sv_n(m,x,s) \quad (37)$$

with initial values $v_0(m,x,s)=m$ and $v_i(m,x,s)=x^i$ for $1\leq i \leq m-1$
and

$$w_{n+m}(m,x,s) = xs^{m-2}w_{n+1}(m,x,s) + s^{m-1}w_n(m,x,s) \quad (38)$$

with initial values $w_0(m,x,s)=m$, $w_i(m,x,s)=0$ for $1\leq i \leq m-2$, and $w_{m-1}(m,x,s)=(m-1)s^{m-2}x.$

Let $c_{m,k}(x,s,z)$ be the characteristic polynomial of (the recurrence of) the sequence $\left(p_k(n,m,x,s)\right)_{n\geq 0}$, so that $c_{m,k}(x,s,E)p_k(n,m,x,s)=0$.

The characteristic polynomial of $v_n(m,x+1,x)$ is

$$c_{m,1}(x+1,x,z) = z^{m-1}(z-1) - x(z^{m-1}-1) = (z-1)(z^{m-1}-xz^{m-2}-xz^{m-3}-\cdots-xz-x) \quad (39)$$

and that of $w_n(m,x+1,(-1)^{m-1}x)$ is

$$c_{m,m-1}(x+1,x,z) = z^m - (x+1)x^{m-2}z - (-1)^{m-1}x^{m-1} = \left(x-(-1)^m z\right)\left(z^{m-1} + xz^{m-2} + \cdots + x^{m-2}z - x^{m-2}\right).$$

It is easily verified that

$$c_{m,m-1}(x+1,x,z) = \frac{z^m}{x} c_{m,1}(x+1,x,(-1)^m \frac{x}{z}). \tag{40}$$

This is in accord with the fact that the right hand side of (10) has the characteristic polynomial $c(z) = (z-1)\left(z-(-1)^m x\right)$ satisfies $c\left((-1)^m \frac{x}{z}\right) = \frac{(-1)^m x}{z^2} c(z).$

These results can be formulated as

**Theorem 4**

*For each $m \geq 2$ we have*

$$p_1(n,m,x,s) = v_n(m,x,s) \tag{41}$$

*and*

$$p_{m-1}(n,m,x,s) = (-1)^m w_n(m,x,(-1)^{m-1}s). \tag{42}$$

**2.4. The case $m = 4$.**

Now we want to consider the case $m = 4$ in more detail.
We already know that

$$p_1(n,4,x,s) = \sum_{j=0}^{\lfloor \frac{n-1}{3} \rfloor} (-1)^j \binom{n-3j}{j} \frac{n}{n-3j} x^{n-4j} s^j$$

and

$$p_3(n,4,x,s) = \sum_{j=\lfloor \frac{2n}{3} \rfloor}^{n-1} (-1)^j \binom{n-j}{3n-4j} \frac{n}{n-j} x^{3n-4j} s^j.$$

The polynomials $p_2(n,4,x,s)$ turn out to satisfy a recursion of order 6 instead of order 4. It is given by

$$p_2(n,4,x,s) = sp_2(n-2,4,x,s) + sx^2 p_2(n-3,4,x,s) + s^2 p_2(n-4,4,x,s) - s^3 p_2(n-6,4,x,s) \tag{43}$$

with initial values
$$p_2(0,4,x,s) = -6, p_2(1,4,x,s) = 0, p_2(2,4,x,s) = -2s, p_2(3,4,x,s) = -3sx^2,$$
$$p_2(4,4,x,s) = -6s^2, p_2(5,4,x,s) = -5s^2x^2.$$

To prove this we observe that we know already the generating functions
$$\sum_{n\geq 0} p_1(n,4,x,s)z^n = \frac{4-3xz}{1-xz+sz^4} \text{ and } \sum_{n\geq 0} p_3(n,4,x,s)z^n = \frac{4-xs^2z^3}{1-xs^2z^3+s^3z^4}.$$

This gives $\sum_{n\geq 0} p_1(n,4,x+1,x)z^n = \frac{4-3(x+1)z}{1-(x+1)z+xz^4} = \frac{1}{1-z} + \frac{3-2xz-xz^2}{1-xz-xz^2-xz^3}$ and

$$\sum_{n\geq 0} p_3(n,4,x+1,x)z^n = \frac{4-(x+1)x^2z^3}{1-(x+1)x^2z^3+x^3z^4} = \frac{1}{1-xz} + \frac{3+2xz+x^2z^2}{1+xz+x^2z^2-x^2z^3}.$$

Furthermore we know that $p_2(n,4,x+1,x) = 1 + x^n - p_1(n,4,x+1,x) - p_3(n,4,x+1,x)$.
This gives the generating function

$$\sum_{n\geq 0} p_2(n,4,x+1,x)z^n = -\frac{3-2xz-xz^2}{1-xz-xz^2-xz^3} - \frac{3+2xz+x^2z^2}{1+xz+x^2z^2-x^2z^3}$$
$$= \frac{-6+4xz^2+3x(x+1)^2z^3+2x^2z^4}{1-xz^2-x(x+1)^2z^3-x^2z^4+x^3z^6}$$

and thus also
$$\sum_{n\geq 0} p_2(n,4,x,s)z^n = \frac{-6+4sz^2+3sx^2z^3+2s^2z^4}{1-sz^2-sx^2z^3-s^2z^4+s^3z^6}.$$

We did not find a simple expression for the polynomials $p_2(n,4,x,s)$ themselves.

### 2.5. The case $m=5$.
For $m=5$ we have found the following characteristic polynomials for $p_k(n,5,x+1,x)$:
$$c_{5,1}(x+1,x,z) = z^5 - (1+x)z^4 + x = (z-1)(z^4 - xz^3 - xz^2 - xz - x),$$

$$c_{5,2}(x+1,x,z) = z^{10} - x(1+x)z^7 - x(1+x)^3z^6 - 2x^2z^5 - x^2(1+x)^2z^4 + x^3(1+x)z^2 + x^4$$
$$= (z^4 - xz^3 - xz^2 - xz - x)(z^6 + xz^5 + x(1+x)z^4 + x^2(1+x)z^3 - x^2(1+x)z^2 + x^3z - x^3),$$

$$c_{5,3}(x+1,x,z) = z^{10} + x(1+x)z^8 - x^2(1+x)^2z^6 + 2x^3z^5 - x^3(1+x)^3z^4 + x^4(1+x)z^3 + x^6$$
$$= (z^4 - xz^3 + x^2z^2 - x^3z - x^3)(z^6 + xz^5 + x(1+x)z^4 + x^2(1+x)z^3 - x^2(1+x)z^2 + x^3z - x^3),$$

$$c_{5,4}(x+1,x,z) = z^5 - (1+x)x^3 z - x^4 = (z+x)(z^4 - xz^3 + x^2 z^2 - x^3 z - x^3).$$

Once found these identities can be verified as above by proving (10).

## 3. Some open problems

Extensive computer experiments suggest the following facts:

$\deg c_{m,k}(x+1,x,z) = \binom{m}{k}$; for $k < \dfrac{m}{2}$ the polynomial $c_{m,k}(x+1,x,z)$ factors in the form

$c_{m,k}(x+1,x,z) = v_{m,k-1} v_{m,k}$ where $v_{m,0} = z-1$; for $k > \dfrac{m}{2}$ the factorization is

$c_{m,k}(x+1,x,z) = v_{m,k} v_{m,k+1}$; for $k = \dfrac{m}{2}$ we get $c_{m,k}(x+1,x,z) = v_{m,k-1} v_{m,k+1}$; $\deg v_{m,k} = \binom{m-1}{k}$;

for $m$ odd we have $v_{m,\lfloor \frac{m}{2} \rfloor} = v_{m,\lceil \frac{m}{2} \rceil}$; $v_{m,m} = z + (-1)^{m-1} x$. Furthermore it seems that as for

$k=0$ and $k=1$ also in general $c_{m,m-k}(x+1,x,z)$ is the monic polynomial proportional to

$c_{m,k}(x+1,x,(-1)^m \dfrac{x}{z})$. So it would suffice to compute $v_{m,k}$ for $k < \dfrac{m}{2}$.

But till now I could not prove these general results.